\documentclass[12pt]{article}
 
\usepackage{amsmath,amsfonts}

\begin{document}
\bibliographystyle{abbrv}

\title{Some results concerning maximum R\'{e}nyi entropy distributions}
\author{Oliver Johnson \and Christophe Vignat}
\date{\today}
\maketitle

\newcommand{\blah}[1]{}
\newtheorem{theorem}{Theorem}[section]
\newtheorem{lemma}[theorem]{Lemma}
\newtheorem{proposition}[theorem]{Proposition}
\newtheorem{corollary}[theorem]{Corollary}
\newtheorem{conjecture}[theorem]{Conjecture}
\newtheorem{definition}[theorem]{Definition}
\newtheorem{example}[theorem]{Example}
\newtheorem{remark}[theorem]{Remark}
\newtheorem{condition}{Condition}
\newtheorem{main}{Theorem}
\newtheorem{assumption}[theorem]{Assumption}

\setlength{\parskip}{\parsep}
\setlength{\parindent}{0pt}
\def \outlineby #1#2#3{\vbox{\hrule\hbox{\vrule\kern #1%
\vbox{\kern #2 #3\kern #2}\kern #1\vrule}\hrule}}%
\def \endbox {\outlineby{4pt}{4pt}{}}%
\newenvironment{proof}
{\noindent{\bf Proof\ }}{{\hfill \endbox
}\par\vskip2\parsep}
\newenvironment{pfof}[2]{\removelastskip\vspace{6pt}\noindent
 {\it Proof  #1.}~\rm#2}{\par\vspace{6pt}}
 
\hfuzz20pt

\newcommand{\Section}[1]{\setcounter{equation}{0} \section{#1}}
\newcommand{\var}{{\rm{Var\;}}}
\newcommand{\cov}{{\rm{Cov\;}}}
\newcommand{\bdiy}{\begin{displaystyle}}
\newcommand{\ediy}{\end{displaystyle}}
\newcommand{\tends}{\rightarrow \infty}
\newcommand{\ep}{{\mathbb {E}}}
\newcommand{\pr}{{\mathbb {P}}}
\newcommand{\re}{{\mathbb {R}}}
\newcommand{\I}{\mathbbm{1}}
\newcommand{\N}{{\mathcal{N}}}
\newcommand{\vcc}[1]{{\mathbf {#1}}}
\newcommand{\vc}[1]{{\bf {#1}}}
\newcommand{\mat}[1]{{\bf {#1}}}
\newcommand{\conpr}{\stackrel{p}{\rightarrow}}
\newcommand{\cond}{\stackrel{d}{\rightarrow}}
\newcommand{\condo}{\stackrel{d^{\circ}}{\rightarrow}} 
\newcommand{\q}{{\cal{Q}}}
\newcommand{\expo}[1]{\exp(#1)}
\newcommand{\tr}{{\rm tr}}
\newcommand{\wt}[1]{ {\widetilde{#1}}}

{\bf English abstract}:
We consider the Student-$t$ and Student-$r$ distributions, which 
maximise R\'{e}nyi entropy under a covariance condition. We show that they have
information-theoretic properties which mirror those of the Gaussian 
distributions, which maximise Shannon entropy under the same condition. We 
introduce a convolution which preserves the R\'{e}nyi maximising
family, and show that
the R\'{e}nyi maximisers are the case of equality in a version of the
Entropy Power Inequality. Further, we show that the R\'{e}nyi maximisers
satisfy a version of the heat equation, motivating the definition of 
a generalized Fisher information.

{\bf French title}: Quelques r\'{e}sultats au sujet des distributions \`{a} 
entropie de R\'{e}nyi maximale.

{\bf French abstract}:
Nous consid\'{e}rons les distributions de types Student-$t$ et Student-$r$ qui
maximisent l'entropie de R\'{e}nyi sous contrainte de covariance. Nous
montrons qu'elles poss\`{e}dent des propri\'{e}t\'{e}s informationnelles similaires 
\`{a}
celles des distributions Gaussiennes, lesquelles maximisent l'entropie de
Shannon sous la m\^{e}me contrainte. Nous montrons que ces distributions sont
stables pour un certain type de convolution et qu'elles saturent une
in\'{e}galit\'{e} de la puissance entropique. De plus nous montrons que les lois 
\`{a}
entropie de R\'{e}nyi maximale v\'{e}rifient une \'{e}quation de la chaleur, ce qui
permet de d\'{e}finir une information de Fisher g\'{e}n\'{e}ralis\'{e}e.

{\bf Keywords}: Entropy Power Inequality, Fisher information, heat equation,
maximum entropy, R\'{e}nyi entropy

{\bf Mathematics Subject Classification}: Primary 94A17; Secondary 60E99

\section{Introduction} \label{sec:intro}

It is natural to ask whether the Shannon entropy of
a $n$-dimensional random vector with density $p$, defined as
$$ H(p) = -\int p(\vc{x}) \log p(\vc{x}) d\vc{x},$$
represents the only possible measure of uncertainty. For example, 
R\'{e}nyi \cite{renyi2} introduces axioms on how we would expect such a 
measure to behave, and shows that these axioms are satisfied by
a more general definition, as follows:
\begin{definition} 
Given a probability density $p$ valued on $\re^n$, for $q \neq 1$ define the 
$q$-R\'{e}nyi entropy to be:
$$ H_{q}(p) = \frac{1}{1-q} \log \left( 
\int p(\vc{x})^{q} d\vc{x}  \right).$$
\end{definition}
Note that by L'H\^{o}pital's rule, since $\frac{d}{dt} a^t = a^t \log_e a$,
\begin{equation} \label{eq:entcont}
\lim_{q \rightarrow 1} H_{q}(p) =
\lim_{q \rightarrow 1} \frac{ - \int p(\vc{x})^{q} \log p(\vc{x}) d\vc{x}}{ 
\int p(\vc{x})^q d\vc{x}} =  H(p). \end{equation}
As Gnedenko and Korolev \cite{gnedenko2} remark, under  a variety of  
natural conditions the  distributions which
maximise  Shannon entropy  are  well-known ones, with interesting properties.  
This paper gives parallels to some of these properties for the R\'{e}nyi 
maximisers.
\begin{enumerate}
\item Under a covariance constraint Shannon entropy  is 
maximised   by  the Gaussian distribution. In Proposition \ref{prop:renmax}
we review the fact that under a covariance constraint R\'{e}nyi entropy is
maximised by Student distributions.
\item The Gaussians have the
appealing property of stability (that is, given $Z_1$ and $Z_2$ 
Gaussians, $Z_1+Z_2$ is also Gaussian). In Definition \ref{def:starconv}, we
introduce the $\star$-convolution, which generalizes the addition operation.
In Lemma \ref{lem:stab}, we extend the stability property by showing that 
if $R_1$ and $R_2$ are R\'{e}nyi maximisers then so is $R_1 \star R_2$.
\item The Entropy Power 
Inequality (see Equation (\ref{eq:epi}) below) shows that the Gaussian
represents the extreme case for how much entropy can change on addition.
Theorem \ref{thm:epi} 
gives the equivalent of an Entropy Power Inequality, with the R\'{e}nyi
maximisers playing an extremal role. 
\item
The Gaussian density satisfies the heat equation, which leads to a 
representation of
Shannon entropy as an integral of Fisher Informations (known as the 
de Bruijn identity).
In Theorem \ref{thm:qheat} we show that the R\'{e}nyi densities satisfy
a generalization of the heat equation, and deduce what quantity must replace 
the Fisher information in general.
\end{enumerate}
First, 
as in Costa, Hero and Vignat \cite{costa2}, we identify the R\'{e}nyi 
maximising densities, which are Student-$t$ and Student-$r$ distributions, and
review some of their properties which we will use later in the paper.
\begin{definition} \label{def:renmax}
For $n/(n+2) < q$ and $q \neq 1$, 
define the $n$-dimensional probability density $g_{q,\mat{C}}$ as
\begin{eqnarray}
g_{q,\mat{C}}(\vc{x}) & = & A_q 
\left( 1- (q-1) \beta \vc{x}^{T}\mat{C}^{-1} \vc{x} \right)_{+}^{\frac{1}{q-1}} 
\label{eq:renmax}
\end{eqnarray}
with \[
\beta= \beta_q = \frac{1}{2q-n\left(1-q\right)},\]
and normalization constants$$
A_{q}= \left\{ \begin{array}{ll}
\left( \Gamma \left( \frac{1}{1-q} \right) (\beta(1-q))^{n/2} \right)/ \left(
\Gamma\left(\frac{1}{1-q}-\frac{n}{2}\right)\pi^{n/2}|\mat{C}|^{\frac{1}{2}} 
\right) 
& \mbox{ if } \frac{n}{n+2}<q<1 \\
\left( \Gamma\left(\frac{q}{q-1}+\frac{n}{2}\right) (\beta(q-1))^{n/2} \right)/
\left( \Gamma\left(\frac{q}{q-1}\right) \pi^{n/2} |\mat{C}|^{\frac{1}{2}} 
\right)
& \mbox{ if } q>1.\\ \end{array} \right. $$
Here $x_+ = \max(x,0)$ denotes the positive part.
We write $\vc{R}_{q,\mat{C}}$
for a random variable with density $g_{q,\mat{C}}$, which 
has mean $\vc{0}$ and covariance $\mat{C}$.
\end{definition}
Notice that if we write $\Omega_{q,\mat{C}}$ for the support of
$g_{q,\mat{C}}$, then for $q > 1$, $\Omega_{q,\mat{C}} = \{ \vc{x}: \vc{x}^T
\mat{C}^{-1} \vc{x} \leq 2q/(q-1) + n \}$, and for $q < 1$, $\Omega_{q,\mat{C}}
= \re^n$. 

Note further that since 
$\lim_{q \rightarrow 1} \Gamma(1/(1-q))(1-q)^{n/2}/\Gamma
(1/(1-q) - n/2) =1$ and
$\lim_{q \rightarrow 1} 
\left( 1- (q-1) \beta \vc{x}^{T}\mat{C}^{-1} \vc{x} 
\right)_{+}^{\frac{1}{q-1}} 
= \exp( - \vc{x}^T \mat{C}^{-1} \vc{x}/2)$, the limit\\
$\lim_{q \rightarrow 1} g_{q,\mat{C}}(\vc{x}) = 
g_{1,\mat{C}}(\vc{x}) = ( (2 \pi)^n |\mat{C}|)^{-1/2} 
\exp( - \vc{x}^T \mat{C}^{-1} \vc{x}/2)$, the Gaussian density.
Throughout this paper, we write $\vc{Z}_{\mat{C}}$ for a $\N(\vc{0}, \mat{C})$
random variable.

We now state the maximum entropy property, as follows.

\begin{proposition} \label{prop:renmax}
Given any $q > n/(n+2)$, and positive definite symmetric 
matrix $\mat{C}$, among all probability
densities $f$ with mean $\vc{0}$ and $\int_{\Omega_{q,\mat{C}}} 
f(\vc{x}) \vc{x} \vc{x}^T d\vc{x} =  \mat{C}$, 
the R\'{e}nyi entropy is uniquely
maximised by $g_{q,\mat{C}}$, that is 
$$H_q(f) \leq H_q(g_{q,\mat{C}}),$$ 
with equality if and only if $f = g_{q,\mat{C}}$ almost everywhere.
\end{proposition}
\begin{proof} See Section \ref{sec:maxent} \end{proof}

Throughout this paper, we write 
$\chi_m$ for a random variable with density 
\begin{equation} \label{eq:chidensity}
f_m(x) = \frac{2^{1-m/2}}{\Gamma(m/2)} x^{m-1} \exp \left(- \frac{x^2}{2}
\right),
\mbox{ for $x > 0$}. \end{equation}
(Strictly speaking, this is only a $\chi$ random variable when the parameter 
$m$ is an integer, but it is simpler to adopt the convention of allowing
non-integer $m$ than to refer to the square root of a $\Gamma(m)$ random 
variable with scale factor 2).

We briefly review
stochastic representations of the R\'{e}nyi maximisers, which we
will use throughout the paper. For the sake of completeness, we present proofs
of these results in Section \ref{sec:stochrep}. Part 1. of Proposition 
\ref{prop:stochrep} follows for example from P.393 of Eaton \cite{eaton},
Part 2. of Proposition \ref{prop:stochrep} is stated in Dunnett
\cite{dunnett}, and Part 3. of this proposition is a multivariate version of
a result stated as long ago as 1915 by Fisher \cite{fisher}.
\begin{proposition} \label{prop:stochrep}
Writing  $\vc{R}_{q,\mat{C}}$   for  a  $n$-dimensional  $q$-R\'{e}nyi
maximiser  with mean  $\vc{0}$ and  covariance $\mat{C}$,  and writing
$\vc{Z}_{\mat{C}}$ for a $\N(\vc{0}, \mat{C})$:
\begin{enumerate} 
\item {\bf Student-$r$.} For any $q > 1$, writing $m = n + 2q/(q-1)$ 
\begin{equation} \vc{R}_{q,\mat{C}} U \sim \vc{Z}_{m \mat{C}},
\label{eq:stochastic_r}\end{equation}
where $U \sim \chi_m$ (independent of
$\vc{R}_{q,\mat{C}}$).
\item {\bf Student-$t$.} For any $n/(n+2) < q <1$, writing 
$m= 2/(1-q) - n > 2$,
\begin{equation} \vc{R}_{q,\mat{C}} \sim \vc{Z}_{(m-2) \mat{C}}/U,
\label{eq:stochastic_t}\end{equation}
where $U \sim \chi_m$ (independent of $\vc{Z}$).
\item {\bf Duality.} 
Given matrix $\mat{D}$, define the map 
$$ \Theta_{\mat{D}}(\vc{x}) = \frac{\vc{x}}{
\sqrt{\vc{x}^T \mat{D}^{-1} \vc{x} + 1}},$$ 
For $q <1$, writing $m = 2/(1-q) - n$, if $\vc{R}_{q,\mat{C}}$ is a 
R\'{e}nyi maximiser,
then $\Theta_{\mat{C}(m-2)}(\vc{R}_{q,\mat{C}}) 
\sim \vc{R}_{p,\mat{C}^*}$, 
where 
$1/(p-1) = 1/(1-q) - n/2-1$ (so $q < 1$ implies that $p > 1$) and
$\mat{C}^* = \mat{C}((m-2)/(m+n))$.
\end{enumerate}
\end{proposition}
\begin{proof} See Section \ref{sec:stochrep}. \end{proof}
Stochastic representations (\ref{eq:stochastic_r}) and (\ref{eq:stochastic_t})
can be used to compute the covariance and entropy of $\vc{R}_{q,\mat{C}}$.
For example, for $q < 1$, since $U \sim \chi_m$, the $ \bdiy
\ep \frac{1}{U^{2}} =\frac{1}{m-2}, \ediy $
 so that $ \bdiy \cov(\vc{R}_{q,\mat{C}}) = \ep \vc{Z}_{(m-2) \vc{C}} 
\vc{Z}^{T}_{(m-2) \vc{C}} \ep \frac{1}{U^{2}}
= (m-2) \mat{C} \frac{1}{m-2}, \ediy$ as claimed.

Similarly for $q < 1$, the Shannon entropy 
$H_{1}\left(\vc{R}_{q,\mat{C}} \right)$ is given by (writing $m = 2/(1-q) - n$)
\begin{eqnarray*}
- \ep \log g_{q,\mat{C}}\left(\vc{R}_{q,\mat{C}} \right) & = &
- \log A_{q} + \frac{m+n}{2} \ep \log\left(1+ \frac{\vc{Z}_{(m-2)\mat{C}}^{T} 
\mat{C}^{-1} 
\vc{Z}_{(m-2) \mat{C}}}{(m-2)U^2} \right)\\
 & = & - \log A_{q} + \frac{m+n}{2}\ep \log\left(1+\frac{\vc{N}^{T} \vc{N}}
{U^2}\right) \\
& = & - \log A_{q} + \frac{m+n}{2} \ep \left(\log \chi^2_{m+n}-\log\chi^2_{m}
\right)
\end{eqnarray*}
where $\vc{N} \sim \N(\vc{0}, \mat{I})$,
and since $\ep \log\chi^2_{m}=\Psi\left(\frac{m}{2}\right)$ where
$\Psi(\cdot)$ is the digamma function, we obtain
\begin{equation} \label{eq:entren}
H_{1}\left( \vc{R}_{q,\mat{C}} \right)=
- \log A_{q} + \frac{1}{1-q}\left(\Psi\left(\frac{1}{1-q}\right)-
\Psi\left(\frac{1}{1-q} - \frac{n}{2}\right)\right). \end{equation}
\begin{remark}
Indeed, the theory of such stochastic representations can be generalized
from the setting of \cite{lutwak} and \cite{lutwak2}
to multivariate maximizers with different powers. That is, given a positive
sequence $\left(p_{1},\dots,p_{n}\right)$, the solution to the problem
\[
\max H_{q}\left(\vc{X}\right)\text{such that } \ep \vert
X_{i}\vert^{p_{i}}=K_{i}
\]
is a random vector $\vc{X}$ with density given by
\[
f\left(\vc{x} \right)\propto\left(1 + \sum_{i=1}^{n}a_{i}\vert
x_{i}\vert^{p_{i}}\right)_{+}^{\frac{1}{q-1}},
\]
where it can be shown that the $a_{i}$ all have the same sign as
$1-q$. Moreover, if $\vc{X}$ is such a maximizer with $q>1,$ then 
for $k = 1, \ldots, n$ random variables
$Z_{k}=U_{k}^{1/p_{k}}X_{k}$ 
are independently power-exponential distributed  with marginal
densities
\[
f\left( z_k \right) = 
\frac{p_k a_k^{\frac{1}{p_k}}}{2 \Gamma (\frac{1}{p_k})}
 \exp \left(a_{k}\vert z_k \vert^{p_{k}} \right), \; \; a_{k} < 0
\]
when $U_{k}$ is $\chi$-distributed with
$m = 2/(q-1)+ 2+\sum_{i=1}^{n} 2/p_{i}$
degrees of freedom and independent of $\vc{X}$.
\end{remark}
\section{$\star$-convolution and relative entropy} \label{sec:epi}
In this section, we introduce a new operation, which we refer to as
the $\star$-convolution. In Lemma \ref{lem:stab} we show that this 
$\star$-convolution preserves the class of R\'{e}nyi entropy maximisers,
and in Theorem \ref{thm:epi} show that it satisfies a version of the entropy
power inequality.

We will say that a distribution is $q$-R\'{e}nyi if it maximises the 
$q$-R\'{e}nyi entropy.
For the sake of simplicity, we write $D(X \| Y) = D_1(f_X \| f_Y)$ for 
the relative entropy
between the two densities $f_X$ and $f_Y$ of random variables $X$ and $Y$.
We define a new distance measure:
\begin{definition} \label{def:newdist}
Given a $n$-dimensional random vector $\vc{T}$ with mean $\vc{0}$ and 
covariance $\mat{C}$, we define its distance from
a $n$-dimensional $q$-R\'{e}nyi maximiser $\vc{R}_{q,\mat{C}}$ (for $q>1$) to be
$$ d(\vc{T} | \vc{R}_{q,\mat{C}}) = D( \vc{T} U \| \vc{Z}),$$
where $U$ is a $\chi_m$ random variable (with $m = n + 2q/(q-1)$ degrees of
freedom)
independent of $\vc{T}$, and $\vc{Z} \sim \N(0, m\mat{C})$. 
\end{definition}
Note that $d$ inherits positive definiteness from $D$ -- that is $d(\vc{T} | 
\vc{R}_{q,\mat{C}})
\geq 0$, with equality if and only if $\vc{T} \sim \vc{R}_{q,\mat{C}}$. 
Note further that Equation (\ref{eq:mult}) below implies that
$$d(\vc{T} | \vc{R}_{q,\mat{C}}) = D( \vc{T} U \| \vc{R}_{q,\mat{C}} U) 
\leq D(\vc{T} \| \vc{R}_{q,\mat{C}}).$$
Motivated by Proposition \ref{prop:stochrep}, we make the following definition:
\begin{definition} \label{def:starconv}
For fixed $q > 1$, given two $n$-dimensional random 
vectors $\vc{S},\vc{T}$,
with covariance matrices  $\mat{C}_{\vc{S}}$ and $\mat{C}_{\vc{T}}$, 
define the 
$\star_q$-convolution (or just $\star$-convolution) of $\vc{S}$ and $\vc{T}$ 
to be the $n$-dimensional random vector
\begin{eqnarray*}
 \vc{S} \star \vc{T} & = & \Theta_{m \mat{C}} \left( \frac{U^{(S)} \vc{S}
+ U^{(T)} \vc{T}}{V} \right) \\
& = &
\frac{(U^{(S)} \vc{S} + U^{(T)} \vc{T})}
{\sqrt{ (U^{(S)} \vc{S} + U^{(T)} \vc{T})^T (m\mat{C})^{-1}
(U^{(S)} \vc{S} + U^{(T)} \vc{T}) + V^2}},
\end{eqnarray*}
where $\mat{C} = \mat{C}_{\vc{S}} + 
\mat{C}_{\vc{T}}$, and $U^{(S)},U^{(T)},V$ are 
independent $\chi$ random variables, where $U^{(S)}$ and $U^{(T)}$ have
$m = n + 2q/(q-1)$ degrees of freedom, and $V$ has $2q/(q-1)$ degrees of 
freedom.
\end{definition}
Again, notice that as $q \rightarrow 1$, $U^{(\cdot)}/(2q/(q-1)) \rightarrow
1$ and $V/(2q/(q-1)) \rightarrow 1$
by the Law of Large Numbers, so $\vc{S} \star \vc{T} \cond \vc{S} + \vc{T}$.
\begin{lemma} \label{lem:stab}
For $q > 1$, if $\vc{S}$ and $\vc{T}$ are $q$-R\'{e}nyi entropy
maximisers with covariances
$\mat{C}_{\vc{S}}$ and $\mat{C}_{\vc{T}}$ then 
$\vc{S} \star \vc{T}$ is also a $q$-R\'{e}nyi entropy maximiser, with 
covariance $\mat{C}_{\vc{S}} + \mat{C}_{\vc{T}}$.
\end{lemma}
\begin{proof}
By Proposition 
\ref{prop:stochrep}.1, writing $m= n + 2q/(q-1)$, we know 
that $U^{(S)} \vc{S}$ and $U^{(T)} \vc{T}$ are $\N(\vc{0},m\mat{C}_{\vc{S}})$ and
$\N(\vc{0},m\mat{C}_{\vc{T}})$ respectively. 
We define $\wt{q}$ by $1/(1- \wt{q}) = 1 + 1/(q-1) + n/2$, and write 
$\wt{m} = 2/(1-\wt{q}) - n = 2q/(q-1) = m-n$. Then random variable
$\vc{W} 
= \sqrt{(\wt{m}-2)/m} (U^{(S)} \vc{S} + U^{(T)} \vc{T})$ 
is $\N( \vc{0}, (\wt{m}-2) \mat{C})$, where $\mat{C} = \mat{C}_{\vc{S}} +
\mat{C}_{\vc{T}}$. 

Then (by Proposition \ref{prop:stochrep}.2) since $V$ has $\wt{m}$ degrees
of freedom, $\vc{W}/V$ is $\wt{q}$-R\'{e}nyi, with covariance $\mat{C}$.
Finally (by Proposition \ref{prop:stochrep}.3), 
$\Theta_{(\wt{m}-2) \mat{C}}(\vc{W}/V)$ is $q^*$-R\'{e}nyi,
where $1/(q^*-1) = 1/(1-q_1) - n/2 - 1 = 1/(q-1)$, so in fact it is 
$q$-R\'{e}nyi with covariance $\mat{C}(\wt{m}-2)/(\wt{m}+n)$. 
Hence, $\vc{S} \star \vc{T}
= \sqrt{m/(\wt{m}-2)} \Theta_{(\wt{m}-2)\mat{C}}(\vc{W}/V)$ is $q$-R\'{e}nyi with
covariance $\mat{C} m/(\wt{m} + n) = \mat{C}$, and the result follows.
\end{proof}
We now give a new ($\star$-convolution) version of the 
classical Entropy Power Inequality, which was first stated by 
Shannon as Theorem 15 of \cite{shannon}, with a `proof' sketched in Appendix 
6. More rigorous proofs appeared in Blachman \cite{blachman} and later
in Dembo, Cover and Thomas \cite{dembo}. The result gives that for
independent $n$-dimensional random vectors $\vc{X}$ and $\vc{Y}$, 
\begin{equation} \label{eq:epi}
\expo{2H(\vc{X+Y})/n} \geq \expo{2H(\vc{X})/n} + \expo{2H(\vc{Y})/n},
\end{equation}
with equality if and only if $\vc{X}$ and $\vc{Y}$ are Gaussian with
proportional covariance matrices.

Writing $\mat{C}_{\vc{X}}$ for the covariance matrix of $\vc{X}$, 
we know that $D( \vc{X} \| \vc{Z}_{\vc{X}}) = (
n \log(2 \pi e) + \log |\mat{C}_\vc{X}|)/2 - H(\vc{X})$, so that the 
Entropy Power Inequality (\ref{eq:epi}) is equivalent to 
\begin{eqnarray} 
\lefteqn{|\mat{C}_{\vc{X}} + \mat{C}_{\vc{Y}}|^{1/n} 
\expo{-2D(\vc{X}+\vc{Y} \| \vc{Z}_{\mat{C}_{\vc{X}} + \mat{C}_{\vc{Y}}})/n} } 
\nonumber \\
& \geq & 
|\mat{C}_{\vc{X}}|^{1/n} \expo{-2D(\vc{X} \| \vc{Z}_{\mat{C}_{\vc{X}}})/n} +
|\mat{C}_{\vc{Y}}|^{1/n} \expo{-2D(\vc{Y} \| \vc{Z}_{\mat{C}_{\vc{Y}}})/n}. 
\label{eq:relepi}
\end{eqnarray}
We give an equivalent of Equation (\ref{eq:relepi}), 
with the $\star$-convolution replacing
the operation of addition.
\begin{theorem} \label{thm:epi} Given $q > 1$,
for independent $n$-dimensional random 
vectors $\vc{S},\vc{T}$ with mean $\vc{0}$ and covariances $\mat{C}_{\vc{S}}$,
$\mat{C}_{\vc{T}}$, 
\begin{eqnarray*}
\lefteqn{
|\mat{C}_{\vc{S}} + \mat{C}_{\vc{T}} |^{1/n} 
\expo{-2d(\vc{S} \star \vc{T} |\vc{R}_{q,\mat{C}_{\vc{S}} + 
\mat{C}_{\vc{T}}})/n}} \\
& \geq &
|\mat{C}_{\vc{S}}|^{1/n} \expo{-2d(\vc{S} | \vc{R}_{q,\mat{C}_{\vc{S}}})/n} +
|\mat{C}_{\vc{T}}|^{1/n} \expo{-2d(\vc{T} | \vc{R}_{q,\mat{C}_{\vc{T}}})/n},
\end{eqnarray*} 
with equality if and only if 
$\vc{S}$ and $\vc{T}$ are $q$-R\'{e}nyi with
proportional covariance matrices.
\end{theorem}
\begin{proof}
By Proposition \ref{prop:mapgood} below
 we know that for $U^{(S)}, U^{(T)}, V,W$ all 
independent and $\chi$-distributed, where $U^{(S)}, U^{(T)}, W$ have 
$m = n + 2q/(q-1)$ degrees of freedom, and $V$ has $2q/(q-1)$ degrees of freedom:
\begin{eqnarray}
\lefteqn{d( \vc{S} \star \vc{T} | \vc{R}_{q,\mat{C}_{\vc{S}}
+ \mat{C}_{\vc{T}}})} \nonumber \\
& = & D( (\vc{S} \star \vc{T}) W \| \vc{Z}_{m (\mat{C}_{\vc{S}}
+ \mat{C}_{\vc{T}})}) \nonumber \\
& =& D \left( \frac{(U^{(S)} \vc{S} + U^{(T)} \vc{T})}
{\sqrt{ (U^{(S)} \vc{S} + U^{(T)} \vc{T})^T \mat{C}^{-1}
(U^{(S)} \vc{S} + U^{(T)} \vc{T}) + V^2}} W \Big\| \vc{Z}_{m 
(\mat{C}_{\vc{S}} + \mat{C}_{\vc{T}})} \right) \nonumber \\
& \leq & D \left( U^{(S)} \vc{S} + U^{(T)} \vc{T} \| \vc{Z}_{m 
(\mat{C}_{\vc{S}} + \mat{C}_{\vc{T}}) } 
\right). \label{eq:dec1}
\end{eqnarray}
We can combine Equations (\ref{eq:relepi}) and 
(\ref{eq:dec1}) to obtain that
\begin{eqnarray*}
\lefteqn{|m \mat{C}_{\vc{S}} + m \mat{C}_{\vc{T}} |^{1/n} 
\expo{-2d(\vc{S} \star \vc{T} |\vc{R}_{q,\mat{C}_{\vc{S}}
+ \mat{C}_{\vc{T}}})/n}} \\
& \geq & |m \mat{C}_{\vc{S}} + m \mat{C}_{\vc{T}} |^{1/n} 
\expo{-2D( U^{(S)} \vc{S} + U^{(T)} \vc{T} \| \vc{Z}_{m (\mat{C}_{\vc{S}} + 
\mat{C}_{\vc{T}})} )/n}  \\
& \geq &
|m \mat{C}_{\vc{S}}|^{1/n} \expo{-2D(U^{(S)} \vc{S} \| 
\vc{Z}_{m \mat{C}_{\vc{S}}})/n}  \\
&  & + |m \mat{C}_{\vc{T}}|^{1/n} \expo{-2D(U^{(T)}\vc{T} \| 
\vc{Z}_{m \mat{C}_{\vc{T}}})/n} 
\\
& = & 
|m \mat{C}_{\vc{S}}|^{1/n} \expo{-2d(\vc{S} | \vc{R}_{q,\mat{C}_{\vc{S}}})/n} +
|m \mat{C}_{\vc{T}}|^{1/n} \expo{-2d(\vc{T} | \vc{R}_{q,\mat{C}_{\vc{T}}})/n},
\end{eqnarray*}
and the result follows. Equality holds in Equation (\ref{eq:dec1}) if 
$U^{(S)} \vc{S}
+ U^{(T)} \vc{T}$ is Gaussian. This, along with proportionality of covariance
matrices, is also the condition for equality in Equation (\ref{eq:relepi}).
\end{proof}
There is a parallel theory for the case $q < 1$, where we define a 
$\circ$-convolution:
\begin{definition} \label{def:starconv2}
For fixed $q$ satisfying $n/(n+2) < q < 1$, given two random vectors
$\vc{S}$ and $\vc{T}$ with covariance matrices $\mat{C}_{\vc{S}}$ and
$\mat{C}_{\vc{T}}$ respectively,
define the $\circ$-convolution by 
\[ \vc{S} \circ \vc{T} = \Theta_{(m-2) (\mat{C}_{\vc{S}} + 
\mat{C}_{\vc{T}})}^{-1} \biggl(
\Theta_{\left(m-2\right) \mat{C}_{\vc{S}}} \left( \vc{S} \right) \star_{\wt{q}}
\Theta_{\left(m-2\right) \mat{C}_{\vc{T}}} \left( \vc{T} \right) \biggr)\]
with $m = 2/(1-q) -n $, where the 
$\star$-convolution is taken with respect to index $\wt{q}$ satisfying 
$1/(\wt{q} - 1) = m/2-1$ and
$$ \Theta_{\mat{D}}^{-1}(\vc{X})=\frac{\vc{X}}
{\sqrt{1- \vc{X}^T \mat{D}^{-1} \vc{X}}}.$$
\end{definition}
This definition satisfies an analogue of Lemma \ref{lem:stab}:
\begin{lemma} \label{lem:stab2}
For $q < 1$, if $\vc{S}$ and $\vc{T}$ are $q$-R\'{e}nyi entropy
maximisers with covariances
$\mat{C}_{\vc{S}}$ and $\mat{C}_{\vc{T}}$ then 
$\vc{S} \circ \vc{T}$ is also a $q$-R\'{e}nyi entropy maximiser, with 
covariance $\mat{C}_{\vc{S}} + \mat{C}_{\vc{T}}$.
\end{lemma}
\begin{proof} By Proposition \ref{prop:stochrep}.3, $\wt{\vc{S}}=\Theta_{\left(m-2\right) \mat{C}_{\vc{S}}}
\left( \vc{S} \right)$
maximises $\tilde{q}$-R\'{e}nyi entropy with $\tilde{q}>1$ such that $1/(\tilde{q}-1)= 1/(1-q)- n/2-1.$

Moreover, the covariance  matrix of $\wt{\vc{S}}$ is $\mat{C}_{\tilde{\vc{S}}}=\frac{m-2}{m+n} \mat{C}_{\vc{S}}$.
The same result holds for $\vc{T}$ and $\mat{C}_{\tilde{\vc{T}}}=\frac{m-2}{m+n} \mat{C}_{\vc{T}}.$
As a consequence of Lemma \ref{lem:stab}, 
$\wt{\vc{S}}\star_{\wt{q}} \wt{\vc{T}}$ is a $\wt{q}$-R\'{e}nyi 
distribution with covariance  
$\wt{\mat{C}}= 
\mat{C}_{\wt{\vc{S}}} + \mat{C}_{\wt{\vc{T}}}.$
 
Since by Proposition \ref{prop:stochrep}.3, $\Theta_{(m-2) \mat{C}}
(\vc{R}_{q,\mat{C}}) = \vc{R}_{\wt{q},\wt{\mat{C}}}$, where 
$\wt{\mat{C}} = (m-2) \mat{C}/(m+n)$,
taking inverse maps, 
$\Theta^{-1}_{(m-2) \mat{C}}( \vc{R}_{\wt{q},\wt{\mat{C}}})
= \vc{R}_{q,\mat{C}}$. Here $\mat{C} = \wt{\mat{C}}(m+n)/(m-2) = 
(\mat{C}_{\wt{\vc{S}}} + \mat{C}_{\wt{\vc{T}}})(m+n)/(m-2) = 
\mat{C}_{\vc{S}} + \mat{C}_{\vc{T}}$, as required.
\end{proof}
\section{$q$-heat equation and $q$-Fisher information} \label{sec:qheat}
In this section, we show that the R\'{e}nyi maximising distributions satisfy
a version of the de Bruijn identity. That is, we can define a Fisher information
quantity, and show in Equation (\ref{eq:debruijn})
that it is the derivative of entropy. First,
we compute the exact constants in a result of Compte and Jou \cite{Compte}.
\begin{theorem} \label{thm:qheat}
For a fixed $\mu$, write $f_{\tau}$ for the density of
a $\vc{R}_{q,\tau^\mu \mat{C}}$ random variable. If 
$ \mu = 2/(2 + n(q-1)/2)$ then 
$f_{\tau}$ satisfies a heat equation of the form \[
K_q \frac{\partial}{\partial\tau}f_{\tau} \left( \vc{x} \right)= 
\sum_{k,l} \mat{C}_{kl}
\frac{\partial^2}{\partial x_k \partial x_l} f_{\tau}^{q}\left(\vc{x} \right) 
\]
with 
$$ K_q = A_q^{q-1}
\frac{2q(2 + n(q-1))}{2q + n\left(q-1\right)}.$$
\end{theorem}
\begin{proof}
By Equation (\ref{eq:renmax}), we know that for a general choice of $\mu$:
$$ f_{\tau}\left( \vc{x} \right)= \frac{A_{q}}{\tau^{n \mu/2}} 
\left(1 - \frac{(q-1) \beta \vc{x}^T \mat{C}^{-1} \vc{x}}{\tau^\mu} 
\right)^{\frac{1}{q-1}}, \mbox{  where }
\beta =\frac{1}{2q - n\left(1-q\right)}. $$
First note that
\begin{equation} \label{eq:tdiff}
\frac{\partial}{\partial\tau}f_{\tau} \left(\vc{x} \right) =  f_{\tau} 
\left( \vc{x} \right) 
\left( -\frac{n \mu}{2 \tau} + \frac{\beta \mu \vc{x}^T \mat{C}^{-1} 
\vc{x}}{\tau^{\mu+1}} 
\left( 1 - \frac{(q-1) \beta 
\vc{x}^T \mat{C}^{-1} \vc{x}}{\tau^{\mu}} \right)^{-1} 
\right).
\end{equation}
Further, for any $k$, writing $\mat{A} = \mat{C}^{-1}$:
\begin{eqnarray*}
\frac{\partial}{\partial x_k} f_{\tau}^{q}\left(\vc{x} \right) & = & 
\frac{A_{q}^{q}}{\tau^{nq \mu/2}} 
\left(1 - \frac{(q-1) \beta \vc{x}^T \mat{C}^{-1} \vc{x}}{\tau^{\mu}} 
\right)^{\frac{1}{q-1}}
\left( \frac{-2 q \beta( \mat{A} \vc{x})_k}{\tau^{\mu}} \right). 
\end{eqnarray*}
Hence, for any $k$, $l$:
\begin{eqnarray*}
\lefteqn{\frac{\partial^2}{\partial x_k \partial x_l} f^{q}_{\tau}
\left(\vc{x} \right)}\\
& = & 
\frac{A_{q}^{q}}{\tau^{nq \mu/2}} 
\left(1 - \frac{(q-1) \beta \vc{x}^T \mat{C}^{-1} \vc{x}}{\tau^{\mu}} 
\right)^{\frac{1}{q-1}}
\left( - \frac{2 q \beta \mat{A}_{kl}}{\tau^{\mu}} \right) \\
& & +  \frac{A_{q}^{q}}{\tau^{nq \mu/2}} 
\left(1 - \frac{(q-1) \beta \vc{x}^T \mat{C}^{-1} \vc{x}}{\tau^{\mu}} 
\right)^{\frac{1}{q-1}-1}
\left( \frac{4 \beta^2 q}{\tau^{2\mu}} 
( \mat{A} \vc{x})_k (\mat{A} \vc{x})_l \right) \\
& = & \frac{A_q^{q-1}}{\tau^{n(q-1) \mu/2}} f_{\tau}\left(\vc{x} \right)
\left( \frac{-2 q \beta \mat{A}_{kl}}{\tau^{\mu}} + 
\left( 
\frac{4 q \beta^2( \mat{A} \vc{x})_k (\mat{A} \vc{x})_l}{\tau^{2 \mu}} \right)
\left(1 - \frac{(q-1) \beta \vc{x}^T \mat{C}^{-1} \vc{x}}{\tau^{\mu}} 
\right)^{-1} \right)
 \end{eqnarray*}
Overall, we deduce that
\begin{eqnarray} \label{eq:partsum}
\lefteqn{ \sum_{k,l} \mat{C}_{kl}
\frac{\partial^2}{\partial x_k \partial x_l} f^{q}_{\tau}\left(\vc{x}\right)} 
\nonumber \\ 
& = & 
\frac{A_q^{q-1}}{\tau^{n(q-1) \mu/2}} f_{\tau}\left(\vc{x} \right)
\left( - \frac{2 q \beta n}{\tau^{\mu}} +  
\frac{4 q \beta^2 \vc{x}^T \mat{C}^{-1} \vc{x} }{\tau^{2 \mu}}
\left(1 - \frac{(q-1) \beta \vc{x}^T \mat{C}^{-1} \vc{x}}{\tau^{\mu}} 
\right)^{-1} \right) \nonumber
\end{eqnarray}
so that equating this with Equation (\ref{eq:tdiff}) we obtain:
$$ K_q = \frac{A_q^{q-1}}{\tau^{n (q -1) \mu/2 + \mu-1}}  
\frac{4 q \beta}{\mu}.$$

Now, we want this to not be a function of $\tau$, so take
$ \mu = 2/(2 + n(q-1))$, and substitute for $\beta$  to obtain
$$ K_q = A_q^{q-1}
\frac{2q(2 + n(q-1))}{2q + n\left(q-1\right)},$$
as claimed.
\end{proof}
Note that the value of the exponent $\mu$ coincides with the one given
by Compte and Jou \cite{Compte}. 
Further, as $\lim_{q \rightarrow 1} A_q^{q-1} = 1$, so that
$ \lim_{q \rightarrow 1} K_q = 2$, as we would expect from the de Bruijn 
identity given in Lemma 2.2 of Johnson and Suhov \cite{johnson3}.

We now evaluate the derivative of the R\'{e}nyi entropy, extending the 
de Bruijn identity:
\begin{eqnarray}
 \frac{\partial}{\partial \tau} H_{q}(f_{\tau}) & = & \frac{1}{1-q} 
\frac{ (q-1) \int f_{\tau}(\vc{x})^{q-1} \frac{\partial}{\partial \tau}
f_{\tau}(\vc{x}) d\vc{x} }{ \int f_{\tau}(\vc{x})^{q} d\vc{x}} \nonumber \\
& = & - \frac{K_q^{-1}}{ \int f_{\tau}(\vc{x})^{q} d\vc{x}} 
\sum_{k,l} \mat{C}_{kl}
\int f_{\tau}(\vc{x})^{q-1} \frac{\partial^2}{\partial x_k
\partial x_l} f^q_{\tau}(\vc{x}) d\vc{x} \nonumber \\
& = & \frac{K_q^{-1}}{ \int f_{\tau}(\vc{x})^{q} d\vc{x}} 
\sum_{k,l} \mat{C}_{kl}
\int \frac{\partial}{\partial x_l} f_{\tau}(\vc{x})^{q-1} 
\frac{\partial}{\partial x_k} f^q_{\tau}(\vc{x}) d\vc{x} \nonumber \\
& = & \frac{K_q^{-1} q(q-1)}{ \int f_{\tau}(\vc{x})^{q} d\vc{x}} 
\sum_{k,l} \mat{C}_{kl} 
\int f_{\tau}(\vc{x})^{2q-3} \frac{\partial}{\partial x_l} f_{\tau}(\vc{x}) 
\frac{\partial}{\partial x_k} f_{\tau}(\vc{x}) d\vc{x} \nonumber \\
& = & K_q^{-1} q(q-1) \tr\left( \mat{C} \vc{J}_q(f_{\tau}) \right),
\label{eq:debruijn}
\end{eqnarray}
where we make the following definitions:
\begin{definition} Given probability density $p$, define the $q$-score function
$$ \vcc{\rho}_{q}(\vc{x}) = \nabla p(\vc{x})/
p(\vc{x})^{2-q},$$ and the $q$-Fisher information matrix to be
$$ \mat{J}_q(p) = \frac{\int p(\vc{x}) \vcc{\rho}_{q}(\vc{x}) 
\vcc{\rho}_{q}^T(\vc{x}) d\vc{x}} { \int p(\vc{x})^{q} d\vc{x}}.$$
\end{definition}
Note that the numerator is the case $p=2$, $\lambda =q$ of the
$(p,\lambda)$ Fisher information introduced
in Equation (7) of \cite{lutwak2}.
We establish a multi-dimensional Cram\'{e}r-Rao inequality:
\begin{proposition}
For the Fisher information $\mat{J}_q$ defined above, given a random variable 
with density $p$ and
covariance $\mat{C}$ then
$$ \mat{J}_q(p) - \frac{ \int p(\vc{x})^q d\vc{x}}{q^2} \mat{C}^{-1}$$
is positive definite, with equality if and only if $p = g_{q,\mat{C}}$
everywhere.
\end{proposition}
\begin{proof}
The key is a Stein-like identity, as usual found using integration by parts, since
\begin{eqnarray*}
\int p(\vc{x}) (\vcc{\rho}_{q}(\vc{x}))_l (\mat{A} \vc{x})_k  d \vc{x} 
& = &
\int \frac{\partial}{\partial x_l} p(\vc{x}) p^{q-1}(\vc{x}) 
(\mat{A} \vc{x})_k d\vc{x}\\
& = & \frac{1}{q} \int \frac{\partial}{\partial x_l}
\left( p^{q}(\vc{x}) (\mat{A} \vc{x})_k \right)  d \vc{x} \\
& = &  - \frac{1}{q} \int p^{q}(\vc{x}) A_{kl} d \vc{x}.
\end{eqnarray*}
This means that for any real $c$, the positive definite matrix 
\begin{eqnarray*}
\lefteqn{
\int p(\vc{x}) (\vcc{\rho}_{q}(\vc{x}) + c \mat{A} \vc{x})
(\vcc{\rho}_{q}(\vc{x}) + c \mat{A} \vc{x})^T d\vc{x} } \\
& = & \int p(\vc{x}) \vcc{\rho}_{q}(\vc{x}) \vcc{\rho}_{q}^T(\vc{x}) d\vc{x} + 
2 \frac{c}{q} \mat{A} \int p^q(\vc{x}) d \vc{x} + c^2 \mat{A}. 
\end{eqnarray*}
So we choose $c = \left( \int p^{q}(\vc{x}) d \vc{x} \right)/q$, and the
result follows. Note that equality holds if and only if $p = g_{q,\mat{C}}$
everywhere, since the R\'{e}nyi maximiser has score function
$\vcc{\rho}(\vc{x}) = A_q^{q-1} (-2 \beta) \mat{A} x$, and 
$\int g_{q,\mat{C}}^q(\vc{x}) d\vc{x}/q = \int g_{q,\mat{C}} A_q^{q-1}
(1 - \beta (q-1) \vc{x}^T \mat{A} \vc{x})/q d \vc{x} = A_q^{q-1} 
(1-\beta (q-1)n))/q = A_q^{q-1} (2 \beta)$.
\end{proof}
Now, we can give the extensivity property for Fisher information defined 
in this way:
\begin{lemma} For a compound system of independent random vectors
$\vc{X}$ and $\vc{Y}$, for $q > 1/2$ the $q$-Fisher information
satisfies:
$$ \mat{J}_q(\vc{X}, \vc{Y}) = \left( \begin{array}{cc}
\alpha_q(\vc{Y}) \mat{J}_q(\vc{X}) & 0 \\
0 & \alpha_q(\vc{X}) \mat{J}_q(\vc{Y}) \\ \end{array} \right),$$
where constant $\alpha_q(\vc{X}) = (\int p^{2q-1}_{\vc{X}}(\vc{x}) d\vc{x} )/
(\int p^q_{\vc{X}}(\vc{x}) d\vc{x} )$ and $\alpha_q(\vc{Y})$ similarly.
\end{lemma}
\begin{proof}
We write $p_{\vc{X},\vc{Y}}(\vc{x},\vc{y}) = p_{\vc{X}}(\vc{x}) 
p_{\vc{Y}}(\vc{y})$,
so that (omitting the arguments for clarity), we can express
$$ \nabla p_{\vc{X},\vc{Y}}  = ( p_{\vc{Y}} \nabla p_{\vc{X}}, p_{\vc{X}}
\nabla p_{\vc{Y}}).$$
Then 
\begin{eqnarray*}
\lefteqn{\iint p_{\vc{X},\vc{Y}}^{2q-3} \nabla p_{\vc{X},\vc{Y}} 
\nabla^T p_{\vc{X},\vc{Y}}  } \\
& = & \left( \begin{array}{cc}
\iint p_{\vc{Y}}^{2q-1} p_{\vc{X}}^{2q-3} \nabla p_{\vc{X}} \nabla^T p_{\vc{X}}
& 0 \\
0 & \iint p_{\vc{X}}^{2q-1} p_{\vc{Y}}^{2q-3} \nabla p_{\vc{Y}} \nabla^T p_{\vc{Y}}
\end{array}
\right) \\
& = & \left( \begin{array}{cc}
\int p_{\vc{Y}}^{2q-1} \int p_{\vc{X}}^q \mat{J}_q(\vc{X})
& 0 \\
0 & \int p_{\vc{X}}^{2q-1} \int p_{\vc{Y}}^{q} \mat{J}_q(\vc{Y})
\end{array}
\right), 
\end{eqnarray*}
since for $q > 1/2$, the off-diagonal term 
\begin{eqnarray*} 
\lefteqn{
\left( \int p_{\vc{X}}^{2q-2} \nabla p_{\vc{X}} \right)
\left( \int p_{\vc{Y}}^{2q-2} \nabla p_{\vc{Y}} \right)} \\
& = & \frac{1}{(2q-1)^2}
\left( \int \nabla p_{\vc{X}}^{2q-1} \right)
\left( \int \nabla p_{\vc{Y}}^{2q-1} \right) =  0,
\end{eqnarray*}
since this is a perfect derivative, and since $p_{\vc{X}}(\vc{x}) \rightarrow 0$
as $\vc{x} \tends$.
The result follows since
$$ \iint p_{\vc{X},\vc{Y}}^{q} =
\left( \int p_{\vc{X}}^{q}  \right) 
\left( \int p_{\vc{Y}}^{q} \right).$$
\end{proof}

\appendix
\section{Proofs} \label{sec:proof}
\subsection{Maximum entropy property} \label{sec:maxent}

In this section we give a proof of Proposition \ref{prop:renmax}, which 
shows that $g_{q,\mat{C}}$ are the R\'{e}nyi entropy maximisers.
The proof uses
Lemma 1 of Lutwak, Yang and Zhang \cite{lutwak2}, which extends the 
classical Gibbs inequality, and
is equivalent to Lemma \ref{lem:gibbs} below.
\begin{definition} \label{label:qrelent}
For $q \neq 1$, given $n$-dimensional probability densities $f$ and $g$, 
define 
the relative $q$-R\'{e}nyi entropy distance from $f$ to $g$ to be
\begin{eqnarray*}
D_q(f \| g) 
& = & \frac{1}{1-q} \log \left( \int g^{q-1}(\vc{x}) f(\vc{x}) d\vc{x} \right) 
+ \frac{1-q}{q} H_q(g) - \frac{1}{q} H_q(f).
\end{eqnarray*}
\end{definition}
For $q =1$, we write $D_1(f \| g) = \int f(\vc{x}) \log( f(\vc{x})/g(\vc{x}))
d \vc{x}$ for the standard relative entropy. We justify this as an extension
by continuity;
as $q \rightarrow 1$, as in (\ref{eq:entcont}), 
$D_q(f \| g) \rightarrow
- \int f(\vc{x}) \log g(\vc{x}) d \vc{x} - H_1(f) = D_1(f \| g)$.
\begin{lemma} \label{lem:gibbs}
For any $q > 0$, and for any probability densities $f$ and $g$,
the relative entropy 
$D_q(f \| g) \geq 0$, with equality if and only if $f = g$
almost everywhere. 
\end{lemma}
\begin{proof} The case $q=1$ is well-known. For $q \neq 1$,
as in Lutwak, Yang and Zhang 
\cite{lutwak2}, the result is a direct application of
H\"{o}lder's inequality to $\exp D_q(f \| g)$.
Although \cite{lutwak2} only strictly speaking considers the 
1-dimensional case, the general case is precisely the same. \end{proof}
As with the Shannon maximisers, we use this Gibbs
inequality Lemma \ref{lem:gibbs} to show that the densities of Definition
\ref{def:renmax} really do maximise the R\'{e}nyi entropy. 

\begin{proof}{\bf of Proposition \ref{prop:renmax}}
Since $f$ and $g_{q,\mat{C}}$ have the same covariance matrix, 
\[
\int_{\Omega_{q,\mat{C}}}
\left(\vc{x}^{T}\mat{C}^{-1}\vc{x}\right)
f\left(\vc{x}\right)d\vc{x} = 
\int_{\Omega_{q,\mat{C}}}
\left(\vc{x}^{T}\mat{C}^{-1}\vc{x}\right)
g_{q,\mat{C}}\left(\vc{x}\right)d\vc{x}.\]
This means that for $q \neq 1$
\begin{eqnarray}
\int_{\Omega_{q,\mat{C}}} 
g_{q,\mat{C}}^{q-1}\left(\vc{x}\right) f\left(\vc{x}\right) d\vc{x} & = &
\int_{\Omega_{q,\mat{C}}} 
A_{q}^{q-1}\left(1- (q-1) \beta \vc{x}^{T}\mat{C}^{-1}\vc{x}\right)
f\left(\vc{x}\right) d\vc{x} \nonumber \\
& = & \int_{\Omega_{q,\mat{C}}}
A_{q}^{q-1}\left(1- (q-1) \beta \vc{x}^{T}\mat{C}^{-1}\vc{x}\right)
 g_{q,\mat{C}}\left(\vc{x}\right) d\vc{x} \nonumber \\
 & = &\int_{\Omega_{q,\mat{C}}}
g_{q,\mat{C}}^{q}(\vc{x}) d\vc{x}. \label{eq:ortho} \end{eqnarray}
For $q=1,$ the equivalent of the orthogonality property Equation (\ref{eq:ortho})
is the well-known fact that \[
\int f(\vc{x}) \log g_{1,\mat{C}}(\vc{x}) d\vc{x} 
= \int g_{1,\mat{C}}(\vc{x}) \log g_{1,\mat{C}}(\vc{x}) d\vc{x}. \]
Using Equation (\ref{eq:ortho}) we simply evaluate 
\begin{eqnarray*}
D_q(f \| g_{q,\mat{C}}) 
& = &  \frac{1}{1-q} \log \left( \int g_{q,\mat{C}}^{q-1}(\vc{x}) f(\vc{x}) 
d\vc{x} \right) 
+ \frac{1-q}{q} H_q( g_{q,\mat{C}}) - \frac{1}{q} H_q(f) \\
& = & \frac{1}{q} \left( H_q(g_{q,\mat{C}}) - H_q(f) \right),
\end{eqnarray*}
and so the result follows by Lemma \ref{lem:gibbs}.
\end{proof}
Note that this is an
alternative proof to that given by Costa, Hero and Vignat \cite{costa2}, 
who introduced a
non-symmetric directed divergence measure
$$
D_{q}\left(f \| g\right)=\mbox{sign}\left(q-1\right)
\int_{\Omega_{q,\mat{C}}}\frac{f^{q}(\vc{x})}{q}+\frac{q-1}{q}g^{q}(\vc{x})-
f(\vc{x})g^{q-1}(\vc{x}) d\vc{x}. $$
The approach of \cite{costa2} is similar to that used by Cover and Thomas
\cite[p.234]{cover} in the Gaussian case. 
The general theory of directed divergence measures is discussed by Csiszar
\cite{csiszar7} and by Ali and Silvey \cite{ali}. 

The paper \cite{lutwak} gives more general results concerning the maximum
entropy property, in a more geometric context.
\subsection{Stochastic Representation} \label{sec:stochrep}
\begin{proof}{\bf of Proposition \ref{prop:stochrep}}

1. By Equation (\ref{eq:chidensity}),
since we take $\beta(q-1) = 1/m$ in Equation (\ref{eq:renmax}),
the density of $\vc{R}_{q,\mat{C}} U$ can be expressed as
 \begin{eqnarray*} 
g(\vc{y}) & = & \frac{2^{1- \frac{m}{2}} A_q}{\Gamma(m/2)} 
\int_{0}^{\infty} \frac{1}{x^n} \left(1-  \frac{\vc{y}^T \mat{C}^{-1} \vc{y}}
{m x^2} \right)^{\frac{1}{q-1}} 
x^{m-1} \exp \left( - \frac{x^2}{2} \right) dx \\
& = & \frac{2^{1-m/2}}{\Gamma(m/2)} A_q 
\exp\left( -\frac{\vc{y}^T \mat{C}^{-1} \vc{y}}{2m} \right) K.
\end{eqnarray*}
Here since $m-n-2 = 2/(q-1)$, taking $u^2 = x^2 - \vc{y}^T \mat{C}^{-1} \vc{y}/m$,
so $u du = x dx$: 
\begin{eqnarray*}
K & = & 
\int_{0}^{\infty} \left(1-  \frac{\vc{y}^T \mat{C}^{-1} \vc{y}}{m x^2}
\right)^{\frac{1}{q-1}} x^{m-n-2} 
\exp \left(- \frac{x^2}{2} + \frac{\vc{y}^T \mat{C}^{-1}
\vc{y}}{2m} \right) x dx \\
& = & \int_{0}^{\infty} u^{\frac{2}{q-1}} 
\exp\left( -\frac{u^2}{2} \right) u du 
= 2^{\frac{1}{q-1}} \Gamma \left( \frac{q}{q-1} \right),
\end{eqnarray*}
 and the result follows,
since the constant 
$$ \frac{2^{1-m/2}}{\Gamma(m/2)} A_q K 
=  \frac{1}{ (2\pi m)^{n/2} |\mat{C}|^{\frac{1}{2}} }.$$
since $1-m/2 + 1/(q-1) = -n/2$ and $\beta(q-1) = 1/m$.

2. In the same way, the density of $\vc{Z}_{(m-2) \mat{C}}/U$ can be 
expressed as
\begin{eqnarray*}
\lefteqn{
\frac{2^{1-\frac{m}{2}}}{\Gamma(m/2)} \int_0^{\infty} 
\frac{x^n}{\sqrt{(2 \pi(m-2))^n |\mat{C}|}}
\exp \left( - \frac{ x^2 \vc{y}^T \mat{C}^{-1} \vc{y} }{2(m-2)} \right)  
x^{m-1} 
\exp \left( - \frac{x^2}{2} \right) dx } \hspace*{1.5cm} \\
& = & \frac{2^{1-\frac{m}{2}-\frac{n}{2}}}{\Gamma(m/2)} 
\sqrt{ \frac{\beta(1-q)}{\pi^n |\mat{C}|}}
\int_0^{\infty} \exp \left( - \frac{ (1+d) x^2}{2} \right)  x^{n+m-1} dx \\
& = & \frac{\Gamma((n+m)/2)}{\Gamma(m/2)} \sqrt{ \frac{\beta(1-q)}{\pi^n
|\mat{C}|}}
(1 + d)^{-\frac{n+m}{2}} \\
& = & \frac{ \Gamma( 1/(1-q))}{\Gamma(1/(1-q) - n/2)} 
\sqrt{ \frac{\beta(1-q)}{\pi^n |\mat{C}|}}
\left(1 + \frac{\vc{y}^T \mat{C}^{-1} \vc{y}}{m-2} \right)^{\frac{1}{q-1}},
\end{eqnarray*}
writing $d = 1 + (\vc{y}^T \mat{C}^{-1} \vc{y})/(m-2)$, and using the facts 
that $(m+n)/2 = 1/(1-q)$ and $1/(m-2) = \beta(1-q)$, the result follows.

3. For this choice of parameters, $\vc{X} = \vc{R}_{q,\mat{C}}$ has density
$ A_q ( 1 + \vc{x}^T \mat{D}^{-1} \vc{x} )^{1/(q-1)}$.

If $\vc{Y} = \Theta_{\mat{D}}(\vc{X})$, we can calculate the Jacobian 
$|\partial \vc{X}|/|\partial \vc{Y}|= (1- \vc{Y}^T \mat{D}^{-1} 
\vc{Y})^{-1-n/2}$. 
Then, the standard change-of-variables relation gives that, since
$1 - \vc{Y}^T \mat{D}^{-1} \vc{Y} = (1 + \vc{X}^T \mat{D}^{-1} \vc{X})^{-1}$,
we know that $\vc{Y}$ has density
\begin{eqnarray*}
 g_{\vc{Y}}(\vc{y}) & = & (1- \vc{y}^T \mat{D}^{-1} \vc{y})^{-1-\frac{n}{2}}
g_{\vc{X}}( \Theta_{\mat{D}}^{-1}(\vc{y})). 
\end{eqnarray*}
Thus, in particular, taking $\vc{X} \sim \vc{R}_{q,\mat{C}}$ and $\mat{D} = 
\mat{C}(m-2)$,
we know that
\begin{eqnarray*}
g_{\vc{Y}}(\vc{y}) & = &  
\left( 1- \vc{y}^T \mat{D}^{-1} \vc{y} \right)^{-1-\frac{n}{2}} 
A_q \left( 1- \vc{y}^T \mat{D}^{-1} \vc{y}\right)^{-\frac{1}{q-1}} \\
& = &  A_q
\left( 1- \vc{y}^T \mat{D}^{-1} \vc{y} \right)^{\frac{1}{p-1}}.
\end{eqnarray*}
Since $p > 1$, we know that $\vc{Y}$ has covariance $\mat{D} \beta_p (p-1)
= \mat{D}/(2p/(p-1) + n)^{-1} =  \mat{C}(m-2)/(m+n)$.

Further $A_q = (\beta_q(1-q))^{n/2} \Gamma \left( \frac{1}{1-q} \right)/ \left(
\Gamma\left(\frac{1}{1-q}-\frac{n}{2}\right) \pi^{n/2}|\mat{C}|^{\frac{1}{2}}
\right) $ \\
$ = (\beta_p(p-1))^{n/2} \Gamma \left( \frac{1}{p-1} + \frac{n}{2} + 1 
\right)/ \left( \Gamma\left(\frac{1}{p-1} + 1\right)\pi^{n/2}| (m-2)/(m+n) 
\mat{C}|^{\frac{1}{2}} \right)  = A_{p},$
as required.
\end{proof}
Note that an alternate, stochastic proof of Equation 
(\ref{eq:stochastic_r}) can be deduced from the polar
factorization property of Student-$r$ vectors (see \cite{barthe3} 
for a detailed study): if $\vc{X}$ is orthogonally invariant and
$\vc{X}=r \vc{U}$ where $\vc{U}$ is uniformly distributed
on the sphere, then $r=\Vert \vc{X}\Vert$ and
$\vc{U}= \vc{X}/\Vert\vc{X}\Vert$
are independent. Since $\vc{R}_{q,\mat{C}}$ is the marginal
of a vector $\vc{U}$ uniformly distributed on the sphere, we
deduce that \[
\vc{R}_{q,\mat{C}}=\frac{\sqrt{m}\mat{C}^{1/2}\vc{Z}}
{\sqrt{\vc{Z}^{T}\vc{Z}+\chi_{m-n}^{2}}}\]
where $\vc{Z}$ is a Gaussian vector, and where random variable
$\sqrt{\vc{Z}^{T}\vc{Z}+\chi_{m-n}^{2}}$ is chi distributed
with $m$ degrees of freedom and independent of $\vc{R}_{q,\mat{C}}$.
Thus, multiplying $\vc{R}_{q,\mat{C}}$ by an independent
chi-distributed
random variable with $m$ degrees of freedom yields a Gaussian vector
with covariance matrix $m\mat{C}$, which is exactly Equation 
(\ref{eq:stochastic_r}).

\subsection{Projection results} 
To prove the Entropy Power Inequality, Theorem \ref{thm:epi}, we prove
a technical result, Proposition \ref{prop:mapgood}. This 
relies on two well-known
results, Lemma \ref{lem:chainrule} and Lemma \ref{lem:proj}. 
Firstly as a consequence of the chain rule for relative 
entropy (see for example Theorem 2.5.3 of Cover and Thomas \cite{cover}):
\begin{lemma} \label{lem:chainrule}
For pairs of random variables $(X,Y)$ and $(U,V)$,
$$ D( (X,Y) \| (U,V)) \geq D( X \| U).$$
Equality holds if and only if for each $x$, the random variables $Y | X = x$ 
and  $V | U =x$ have the same distribution. In particular if $(X,Y)$
and $(U,V)$ are independent pairs, equality holds if and only if  
$Y$ and $V$ have the same distribution.
\end{lemma}
Secondly, we recall a projection identity, first stated as Corollary 4.1 of 
\cite{kullback2}:
\begin{lemma} \label{lem:proj}
For random vectors $\vc{X}$ and $\vc{Y}$, and for any
invertible function $\Phi$:
$$ D( \Phi(\vc{X}) \| \Phi(\vc{Y})) = D(\vc{X} \| \vc{Y}).$$ 
\end{lemma}
\begin{proposition} \label{prop:mapgood}
For a $n$-dimensional random vector $\vc{M}$, take $N \sim \chi_{2q/(q-1)}$
and $U \sim \chi_{2q/(q-1)+n}$, where $(\vc{M}, N, U)$
are independent:
$$ D \left( \frac{\vc{M}}{\sqrt{\vc{M}^T \mat{C}^{-1} \vc{M}+N^2}} 
U \Big\| \vc{Z}_{\mat{C}} \right) \leq D( \vc{M} \| \vc{Z}_{\mat{C}}),$$
where  equality holds if $\vc{M}$ is $\N(\vc{0},\mat{C})$.
\end{proposition}
\begin{proof} 
By combining Lemmas \ref{lem:chainrule} and \ref{lem:proj}, if 
random variables $Q$ and $S$ have the
same distribution and $(\vc{P},Q)$ and $(\vc{R},S)$ each form independent pairs
then
\begin{equation} \label{eq:mult} 
D( \vc{P} Q \| \vc{R} S) \leq D( (\vc{P} Q,Q) \| (\vc{R} S, S))
= D( (\vc{P},Q) \| (\vc{R},S)) = D(\vc{P} \| \vc{R}). \end{equation}
Now, we define $Y \sim \chi_{2q/(q-1)}$ and 
$V \sim \chi_{2q/(q-1)+n}$, both independent of $\vc{Z}_{\mat{C}}$, 
so that $U$ and $V$ have the same distribution, as do $N$ and $Y$.
The LHS of the proposition becomes:
\begin{eqnarray}
\lefteqn{ D \left( \frac{\vc{M}}{\sqrt{\vc{M}^T \mat{C}^{-1} \vc{M}+N^2 }} 
U \Big\| 
\frac{\vc{Z}_{\mat{C}}}{\sqrt{\vc{Z}^T_{\mat{C}} \mat{C}^{-1} \vc{Z}_{\mat{C}} + Y^2}} 
V \right)}  
\nonumber \\
& \leq &  D \left( \frac{\vc{M}}{\sqrt{\vc{M}^T \mat{C}^{-1} \vc{M}+N^2}} 
 \Big\| \frac{\vc{Z}_{\mat{C}}}{\sqrt{\vc{Z}^T_{\mat{C}} \mat{C}^{-1} 
\vc{Z}_{\mat{C}} + Y^2}} \right) \label{eq:step2} \\
& = &  D \left( \Theta_{\mat{C}}(\vc{M}/N) \| \Theta_{\mat{C}}
(\vc{Z}_{\mat{C}}/Y) \right) \nonumber \\
& = &  D \left( \vc{M}/N \| \vc{Z}_{\mat{C}}/Y \right) 
\label{eq:step3} \\
& \leq &  D \left( \vc{M} \| \vc{Z}_{\mat{C}} \right) \label{eq:step5},
\end{eqnarray}
and the result follows.
Here Equation (\ref{eq:step2}) follows by Equation (\ref{eq:mult}), 
Equation (\ref{eq:step3}) follows by 
Lemma \ref{lem:proj} and  Equation
(\ref{eq:step5}) again follows by Equation (\ref{eq:mult}).
\end{proof}

\section*{Acknowledgment}
This work was done during a visit by CV to OTJ at the Statistical
Laboratory, University of Cambridge, in April 2005. CV would like to
thank OTJ and the Statistical Laboratory for their hospitality.

\section*{Addresses}

Oliver Johnson, Statistical Laboratory, DPMMS,
Centre for Mathematical
Sciences, University of Cambridge, Wilberforce Road, Cambridge CB3 0WB, UK. 
Fax: +44 1223 337956. Email: {\tt{otj1000@cam.ac.uk}}.

Christophe Vignat
L.I.S., 961 rue de la Houille Blanche, 38402 St. Martin d'H\`{e}res cedex,
France. Email: {\tt{vignat@univ-mlv.fr}}.

\end{document}